\documentclass[a4paper,12pt]{amsart}
\usepackage{ifthen,verbatim}
\usepackage{mathrsfs}

\usepackage{epsfig}
\usepackage{pstricks}
\usepackage{pst-node}
\usepackage{pst-tree}
\usepackage{pst-plot}
\usepackage{pst-text}
\usepackage{multido}
\usepackage{subfigure}

\definecolor{hylightgray}{rgb}{0.9,0.9,0.9}

\numberwithin{equation}{section}
\nonstopmode
\newtheorem{thm}{Theorem}[section]
\newtheorem{cor}[thm]{Corollary}
\newtheorem{lem}[thm]{Lemma}
\newtheorem{prop}[thm]{Proposition}

\theoremstyle{definition}

\newtheorem{example}[thm]{Example}

\newenvironment{rem}{%
\bigskip
\noindent
\textsl{{\sl Remark. }}}{\bigskip}

\newenvironment{pf}[1][]{%
 \vskip 3mm
 \noindent
 \ifthenelse{\equal{#1}{}}%
  {{\slshape Proof. }}%
  {{\slshape #1.} }%
 }%
{\qed\bigskip}

\newcounter{alphabet}
\newcounter{tmp}

\newcommand{\CR}{{\operatorname{cr}}}
\newcommand{\C}{{\mathbb C}}
\newcommand{\D}{{\mathbb D}}

\newcommand{\mQ}{{\widehat Q}}
\newcommand{\R}{{\mathbb R}}

\newcommand{\bD}{{\overline{\mathbb D}}}
\newcommand{\bDs}{{\overline{\mathbb D}\,^*}}
\newcommand{\sphere}{{\widehat{\mathbb C}}}

\renewcommand{\Re}{{\operatorname{Re}\,}}

\newcommand{\id}{{\operatorname{id}}}

\newcommand{\card}{{\operatorname{card}}}

\renewcommand{\mod}{{\operatorname{mod}\,}}
\newcommand{\inv}{^{-1}}

\newcommand{\arccosh}{{\operatorname{arccosh}\,}}

\newcommand{\deltat}{{\hat\delta}}
\newcommand{\Vt}{{\hat V}}
\newcommand{\jt}{{\hat j}}
\newcommand{\qt}{{\hat q}}

\renewcommand{\arg}{\,{\operatorname{arg}\,}}

\newcommand{\sys}{{\operatorname{sys}\,}}

\newcommand{\aand}{{\quad\text{and}\quad}}

\newcounter{minutes}
\divide\time by 60
\newcounter{hours}
\multiply\time by 60
\addtocounter{minutes}{-\time}

\begin{document}
\bibliographystyle{amsplain}
\title{
On the hyperbolic distance of $n$-times punctured spheres
}

\def\thefootnote{}
\footnotetext{
\texttt{\tiny File:~\jobname .tex,
          printed: \number\year-\number\month-\number\day,
          \thehours.\ifnum\theminutes<10{0}\fi\theminutes}
}
\makeatletter\def\thefootnote{\@arabic\c@footnote}\makeatother

\author[T. Sugawa]{Toshiyuki Sugawa}
\address{Graduate School of Information Sciences,
Tohoku University, Aoba-ku, Sendai 980-8579, Japan}
\email{sugawa@math.is.tohoku.ac.jp}
\author[M. Vuorinen]{Matti Vuorinen}
\address{Department of Mathematics and Statistics, University of Turku, FI-20014 Turku, Finland}
\email{vuorinen@utu.fi}
\author[T. Zhang]{Tanran Zhang}
\address{Department of Mathematics, Soochow University, No.1 Shizi Street, Suzhou 215006, China}
\email{trzhang@suda.edu.cn}
\begin{abstract}
The length of the shortest closed geodesic in a hyperbolic surface $X$ is called the systole of $X.$
When $X$ is an $n$-times punctured sphere $\sphere\setminus A$ where $A \subset \sphere$ is a finite
set of cardinality $n\ge4,$ we define a quantity $Q(A)$ in terms of cross ratios of quadruples in $A$
so that $Q(A)$ is quantitatively comparable with the systole of $X.$
We next propose a method to construct a distance function $d_X$ on a punctured sphere $X$
which is Lipschitz equivalent to the hyperbolic distance $h_X$ on $X.$
In particular, when the construction is based on a modified quasihyperbolic metric,
$d_X$ is Lipschitz equivalent to $h_X$ with Lipschitz constant depending only on $Q(A).$
\end{abstract}
\thanks{
The authors were supported in part by JSPS Grant-in-Aid for
Scientific Research (B) 22340025.
}
\maketitle
\section{Introduction}
One of the fundamental research themes of geometric function theory is the study of the way
in which distances between points are transformed or distorted under analytic or meromorphic functions
$f: \Omega_1 \to \Omega_2\,,$ where $\Omega_1$ and $\Omega_2$ are subdomains of the Riemann sphere
$\sphere={\mathbb{C}} \cup \{ \infty \}\,.$
For this purpose distances are measured in terms of metrics. Several metrics are
commonly used: the euclidean, chordal and hyperbolic metrics are some examples.
The numeric value of the hyperbolic distance  $h_\Omega(z,w)$ of two points $z,w$
in a plane domain $\Omega \subset \C$ with $\card(\sphere \setminus \Omega) \ge 3$
does not merely reflect ``the distance between points" but also
``the location of the points with respect to boundary" and ``the metric structure of the boundary".
These facts, along with its conformal invariance, make the hyperbolic metric a powerful
tool in geometric function theory \cite{Ahlfors:conf}. Often it is more natural to use
the hyperbolic distance $h_\Omega$ than to use the euclidean or chordal metric.
The study of the hyperbolic metric is a part
of the wider context where metrics and their applications to geometric function
theory have a key role.  For instance, in the study of quasiconformal maps of subdomains of $\mathbb{R}^n, n\ge2\,,$
several new metrics have been introduced during the past few
decades, with the purpose of generalizing hyperbolic metric to subdomains of $\mathbb{R}^n, n\ge3\,.$
One of these ``hyperbolic-type" metrics is
 the quasihyperbolic metric, which has become a popular tool applied by many authors \cite{GH:qd}, \cite{Heinonen:metric}, \cite{Vai90} etc.

The hyperbolic metric of a general plane domain $\Omega \subset \C$ with
$\card(\sphere \setminus \Omega) \ge 3$ is defined in three steps.
First, in the case of the unit disk $ \Omega = \D$, one can give an explicit formula for the hyperbolic distance \cite{Beardon:disc}.
Second, for a simply connected domain, the definition reduces to the case of the unit disk by use of the Riemann mapping theorem.
Third, in the case of a general domain  $\Omega \subset \C\,,$ the hyperbolic metric is defined as a Riemannian metric
with the density function $\lambda_\Omega=\lambda_\Omega(z)|dz|$  given in terms of
the universal covering mapping (see below or \cite{KL:hg}).
Distances between points are then obtained by minimising the line integrals over the paths, the hyperbolic lengths of the paths,
connecting the points. The reader should notice that in the literature, the word ``metric" refers often to the
density function of the Riemannian metric rather than to the distance function ``metric" in the sense of topology.

Of these three cases, the extremal cases are the first case, when the domain $\Omega$ is the unit disk $\D\,,$
and the third case, when the domain is the
trice punctured sphere $\Omega=\C_{0,1}=\C \setminus\{0,1\}$,
or a domain conformally equivalent to one of these.
The unit disk case is well-known \cite{Beardon:disc} and the density of
the metric is $1/(1-|z|^2)\,.$
The case of  $\C_{0,1}$ is much more difficult and
leads to classical special functions such as elliptic functions and their representation
as infinite products, see \cite{Agard68, Hempel79, KL:hg},
and the density $\lambda(z)=\lambda_{\C_{0,1}}(z)$ of the metric satisfies
\begin{equation} \label{eq:hd}
\lambda(z)=\frac{1}{2|z| |\log |z||}+O(1) \quad\text{as}~ z\to0\,, \quad \lambda(-1) = \frac{2 \pi^2}{\Gamma(1/4)^4} \approx 0.114237\,.
\end{equation}
Functional identities for classical special functions also yield identities and numerical approximations for  $\lambda(z)$ \cite{SV01}.
The case of a general domain lies between these extreme cases which correspond to steps one and three of the above definition.
In the case of the unit disk, we have not only an explicit distance formula, but also a rich structure of the hyperbolic geometry which
includes counterparts of some results of the euclidean geometry  \cite{Beardon:disc}. In the case of a general simply connected domain, we loose
all this because the Riemann mapping function is not explicitly known, but what remains are good upper and lower bounds for the
density \cite[(3.2.1)]{GH:qd}, and if the domain is not only simply connected but also uniform, then the hyperbolic distance has upper and lower
bounds in terms of the distance ratio metric or the $j$-distance,  \cite{GH:qd} .
In the case of a general plane domain neither  good upper nor lower bounds for the
hyperbolic distance exist, only there are upper and lower bounds for the density function of the metric.
These bounds depend on the geometry of the domain and on the metric structure of its boundary.
Several particular cases, such as rectangles, convex domains, quasidisks, domains with isolated boundary points have been studied
in the literature and the  influence of thickness properties of the boundary such as uniform perfectness  on the hyperbolic metric have been
analysed \cite{Bea01, Bea03, BP78, JV96, MM12, SugawaUP}, \cite[pp.~245--257]{KL:hg}. In very few cases,
the density function of the hyperbolic metric is known (annulus, punctured disk, strip domain \cite[pp.~133--135]{KL:hg}).

Many authors have studied these topics which bring together extremal problems of geometric function theory,
classical hyperbolic geometry, special classes of domains,
metric structure conditions of sets, and special functions
\cite{Ahlfors:conf, Agard68, BP78, Bet08, GL01, Hempel79, Hayman:sub2, HMM08, JV96, kmm, KR:cm, KR16, Qiu98, SV01, SV05}.
Mostly, the density function of the distance is studied.
In spite of all this research, finding information about the hyperbolic distance  $h_\Omega(z_1,z_2)$
induced by $\lambda_\Omega(z)|dz|$ in a general domain $\Omega$ seems to be very difficult.
If good explicit bounds were known, one could for instance find new distortion
theorems for analytic functions which are contraction mappings between hyperbolic spaces.
For instance, Landau's and Schottky's theorems for analytic functions are closely connected
with such estimates \cite[pp.~707--708]{Hayman:sub2}, \cite{GU16}, \cite{Hempel79}, \cite[p.~210]{AVV:conf}.
S. Rickman \cite{Rickman84} proved a counterpart of Schottky's theorem for quasiregular mappings
$f: \mathbb{B}^n \to  \mathbb{R}^n \setminus A, n\ge 2\,,$ where $A \subset  \mathbb{R}^n $
is a finite set containing sufficiently many points and where $ \mathbb{B}^n$ is the unit ball of $ \mathbb{R}^n\,.$
For that purpose he constructed a metric, which for $n=2$ is similar to the hyperbolic metric.

In the present paper, we will concentrate on the particular case when $\Omega$ is a
punctured sphere of the form $X=\sphere\setminus A$ for a finite subset $A$ of $\sphere$
with $n=\card(A)\ge 3.$ The aim of this paper is to introduce a comparison function,
a new metric (in the sense of topology) $d_X,$ whose values can be easily computed,
with comparison bounds explicitly controlled by the geometry of $X$ or, equivalently, of $A.$
It is not difficult to observe that the bounds cannot be taken uniformly when $n\ge4.$
Therefore, we have to pay attention to hyperbolic geometry of $X$ in an appropriate way.
To measure ``thickness" of $X,$ it is customary to consider the systole $\sys(X)$ of $X,$
which is defined to be the minimal length of closed hyperbolic geodesics in $X.$
For instance, it is known that $X$ tends to the boundary in the moduli space $\mathcal{M}_{0,n}$
precisely when $\sys(X)\to0.$
However, by the reason described above, it is not easy to compute the systole.

In Section 2 we will introduce a quantity $Q(A)$ which can be computed by means of the cross ratios of quadruples
of points in $A$ in a certain way, see \eqref{QA}.
The main theorem in this section is Theorem \ref{thm:hyp} which roughly states that
$Q(A)$ behaves like $\pi^2/\sys(X)$ when $\sys(X)$ is small.
We believe that this connection might be helpful in other aspects.

During the past few decades, several kinds of hyperbolic-type distances have been proposed such as the quasihyperbolic distance.
However, near the punctures their behaviour is usually very different from that of the hyperbolic distance.
In Section 3, we give a general method to construct a new distance $d_X$ from such a hyperbolic-type distance, say $\mu,$
in such a way that the resulting distance is Lipschitz equivalent to the hyperbolic distance.
We also give a concrete bound for these Lipschitz constants in Theorem \ref{thm:general}.
This construction was first proposed in \cite{SZ17} when $\mu$ is the euclidean distance.

Section 4 is devoted to some applications of this theorem.
As one can easily guess, the euclidean distance is not suitable to construct a nearly hyperbolic distance.
Typical hyperbolic-type distances are the so-called $j$-distance ($j$-metric) and the quasihyperbolic distance.
However, these metrics do not satisfy the assumptions of Theorem \ref{thm:general}.
Thus, we will modify the construction of those metrics and obtain Lipschitz equivalent metrics which
satisfy those assumptions.
In particular, the resulting distance based on the modified quasihyperbolic distance is
Lipschitz equivalent to the hyperbolic distance with bounds depending only on $Q(A)$
(see Theorem \ref{thm:q} below).

\section{Hyperbolic geometry and cross ratio}

We introduce here some basic notation and terminology. We denote by $\D(a,r)$ the open disk $|z-a|<r$ in the complex plane $\C$
for $a\in\C$ and $r>0$ and by $\bD(a,r)$ its closure, namely, the
closed disk $|z-a|\le r.$
Also, let $\D^*(a,r)=\D(a,r)\setminus\{a\}$ and $\bDs(a,r)=\bD(a,r)\setminus\{a\}.$

Let  $\Omega \subset \C$  be a domain with $\card(\sphere \setminus \Omega) \ge 3$ and
with the hyperbolic metric $\lambda_\Omega =\lambda_\Omega(z)|dz|$
of constant Gaussian curvature $-4.$
The hyperbolic length of a piecewise smooth curve $\gamma$ in $\Omega$ is defined to be
$$
\ell_\Omega(\gamma)=\int_\gamma \lambda_\Omega(z)|dz|\,.
$$
We denote by $h_\Omega(z,w)$ the hyperbolic distance on $\Omega$ induced by
$\lambda_\Omega$ and defined for  $z,w \in \Omega$ by
$$
h_\Omega(z,w) = \inf \ell_\Omega(\gamma)
$$
where the infimum is taken over all curves joining $z$ and $w$ in $\Omega.$
It is known that the infimum is attained by a hyperbolic geodesic.
Here, a curve $\gamma$ is called a {\it hyperbolic geodesic}
if its hyperbolic length is minimal among the homotopy class of $\gamma$ keeping the endpoints fixed.
R. Nevanlinna's {\it principle of the hyperbolic metric} \cite[p.~49]{Nev:Ein}, \cite[p.~683]{Hayman:sub2}
implies that if $\Omega_1, \Omega_2$ are domains in $\C$ with at least two finite boundary points and
$f: \Omega_1 \to \Omega_2$ is an analytic function, then
$$
h_{\Omega_2}(f(z),f(w)) \le h_{\Omega_1}(z,w)
$$
for $z,w \in \Omega_1\,,$ i.e. analytic functions are contractions in hyperbolic metric. This contraction property may be regarded
as a very general form of the Schwarz lemma, which is the special case
$\Omega_1 = \Omega_2= \D$. Applying
this result to the inclusion mapping we see that if $\Omega_1 \subset \Omega_2$ and $z,w \in \Omega_1\,,$ then
$$
 h_{\Omega_1}(z,w) \ge  h_{\Omega_2}(z,w) \,.
$$
Hence the hyperbolic metric is monotone with respect to the domain inclusion.

As pointed out in the introduction, the density $\lambda_\Omega(z)$ of the hyperbolic metric
 is not known for general plane domains  $\Omega \subset \C\,.$
 Sometimes it is useful to study another metric with
an explicit density function and use it as a comparison
function for the hyperbolic metric. For such a purpose
we introduce the quasihyperbolic distance $q_\Omega(z_1,z_2)\,.$ The quasihyperbolic distance
is defined for $\Omega\subset\C$ in the same way as the hyperbolic metric through the
length minimizing properties of curves, but
by replacing the hyperbolic density by $1/\delta_\Omega(z),$ where
$$
\delta_\Omega(z)=\min_{a\in\partial\Omega}|z-a|\,.
$$
In the particular case when $\Omega=\C^*=\C\setminus\{0\},$ Martin and Osgood \cite{MO86}
gave an explicit formula of $q(z_1,z_2)=q_{\C^*}(z_1,z_2)$ as follows:
\begin{equation}\label{eq:MO}
q(z_1,z_2)=\sqrt{(\log|z_1|-\log|z_2|)^2
+\Big(2\arcsin\tfrac12\big|\tfrac{z_1}{|z_1|}-\tfrac{z_2}{|z_2|}\big|\Big)^2},\quad
z_1,z_2\in\C^*.
\end{equation}
Note that, if the origin is an isolated boundary point of a hyperbolic domain $\Omega_1\,,$  then for a fixed point
$z_0 \in \Omega_1\,, $ $q_{\Omega_1}(z,z_0)=\log(1/|z|)+O(1)$ as $z\to0$ whereas
$h_{\Omega_1}(z,z_0)=\frac12\log(\log(1/|z|))+O(1)$ as $z\to0$ (see \eqref{eq:hd}).
We also note that by \eqref{eq:MO} $q(tz_1,tz_2)=q(z_1,z_2)$ for $t\in\C^*.$

It is well-known that an $n$-times punctured sphere $X=\sphere\setminus A,$
where $A=\{a_1,\dots,a_n\},$ is  hyperbolic if and only if $n\ge 3.$
In this paper, we will mainly consider hyperbolic punctured spheres,
in other words, the $n$-times punctured sphere
$X=\sphere\setminus A$ with $3\le n=\card(A)<+\infty.$
For convenience, we normalise $X$ so that $0,1,\infty\in A$
in the rest of the paper.
In particular, $X\subset\C.$

We want to have a useful quantity
which can be computed easily by the set of punctures $A$ of $X.$
First let
$$
\mQ(A)=\max_{a_j\in A} \log^+|\CR(a_1,a_2,a_3,a_4)|
$$
for a finite subset $A$ of $\sphere$ with $\card (A)\ge4,$
where $\log^+x=\max\{\log x,0\},$
the maximum is taken over all the quadruples of distinct points in $A$ and
$\CR(a_1,a_2,a_3,a_4)$ stands for the cross ratio; that is,
$$
\CR(a_1,a_2,a_3,a_4)=\frac{(a_1-a_3)(a_2-a_4)}{(a_1-a_2)(a_3-a_4)}
$$
if none of $a_j$ is $\infty,$ and otherwise it is defined as a suitable limit.
Note the identities
\begin{equation}\label{eq:cr}
\CR(a_1,a_2,a_3,a_4)
=\CR(a_2,a_1,a_4,a_3)
=\CR(a_3,a_4,a_1,a_2)
=\CR(a_4,a_3,a_2,a_1)
\end{equation}
and
\begin{equation}\label{eq:cr2}
\CR(a_2,a_1,a_3,a_4)
=1-\CR(a_1,a_2,a_3,a_4).
\end{equation}
We set $\mQ(A)=0$ when $\card (A)<4.$
Let $\sys(X)$ be the infimum of the lengths of closed hyperbolic geodesics in $X.$
This quantity is known as the {\it systole} of the surface $X$ and it gives
a convenient way to measure ``thickness" of the surface $X$
(see \cite{SugawaK, Sugawa10MZ}).
Note that $\sys(X)>0$ for the $n$-times punctured sphere, $n\ge3.$
Schmutz \cite[Theorem 14]{Sch94} (see also \cite[Theorem 7.2]{Adams98}) showed that
\begin{equation}\label{eq:Sch}
\sys(X)\le 2\arccosh(3-6/n)<2\arccosh 3=4\log(1+\sqrt2)
\end{equation}
for an $n$-times punctured sphere $X$ with $n\ge4.$
We note that $\sys(X)$ is not comparable with $\mQ(A)$ in a uniform way
(see Example \ref{ex:2} given below).

We thus consider a modified quantity as described in the following.
A pair of subsets $A_1, A_2$ of $A$ will be called an {\it admissible partition}
of $A$ if each of $A_1$ and $A_2$ contains at least two points and
if $A$ is the disjoint union of $A_1$ and $A_2.$
We write $A_1\sqcup A_2=A$ for it.
For an admissible partition $A=A_1\sqcup A_2,$ we put
$$
P(A_1,A_2)=\min_{a_1,a_1'\in A_1, a_2, a_2'\in A_2}\log^+|\CR(a_1,a_1',a_2,a_2')|.
$$
Note that $P(A_1,A_2)=P(A_2,A_1)$ thanks to \eqref{eq:cr}.
Finally, we set
\begin{equation}\label{QA}
Q(A)=\max_{A_1\sqcup A_2=A}P(A_1,A_2)
\end{equation}
for $A\subset\sphere$ with $4 \le \card (A)< \infty\,,$ where the maximum is taken over all the
admissible partitions of $A.$
Again, we set $Q(A)=0$ if $A$ consists of at most three points.
Obviously,
$$
0\le Q(A)\le \mQ(A).
$$
Note that both $Q(A)$ and  $\mQ(A)$ are invariant under M\"obius transformations.

A doubly connected domain $R$ in $\sphere$ is called an {\it annulus}.
An  annulus
whose boundary consists of two circles (including lines) is called a {\it circular} annulus.
Furthermore, if $R$ is bounded in $\C$ and if $\partial R$ consists of two concentric circles (centered at $a$),
it is called a {\it round} annulus (centered at $a$).
We recall that the modulus $\mod R$ of an annulus
$R$ is defined to be $\log(r_2/r_1)$ when $R$ is conformally equivalent to
the round annulus $\{z\in\C: r_1<|z|<r_2\}.$
We say that an annulus $R$ separates  a set $A_1$ from a set $A_2$
if $A_1$ and $A_2$ are contained in different connected components of $\sphere\setminus R.$
We can understand the geometric meaning of the quantity $P(A_1,A_2)$ in terms of
separating circular annuli as follows.

\begin{lem}\label{lem:mod}
Let $A_1$ and $A_2$ be disjoint finite subsets of $\sphere$ with $\card A_j\ge2$ for $j=1,2.$
If $P(A_1,A_2)>0,$ there exists a circular annulus $R$ separating
$A_1$ from $A_2$ with $\mod R=P(A_1,A_2).$
Conversely, if a circular annulus $R$ separates $A_1$ from $A_2,$ then
$$
P(A_1,A_2)\ge 2\log\sinh\left(\frac{\mod R}2\right)\ge \mod R-2\log(1+\sqrt2).
$$
\end{lem}

\begin{pf}
We first show the former part.
By the M\"obius invariance of the quantities involved, we may assume that
$0\in A_1$ and $\infty\in A_2.$
For $a_1\in A_1, a_2\in A_2$ with $a_1\ne0, a_2\ne\infty,$ we have
$|\CR(a_1,0,\infty,a_2)|=|a_2/a_1|.$
Since $\log^+|\CR(a_1,0,\infty,a_2)|\ge P(A_1,A_2)>0,$ we obtain
$|a_2/a_1|=|\CR(a_1,0,\infty,a_2)|\ge e^{P(A_1,A_2)}>1.$
Set $r_1=\max\{|a_1|: a_1\in A_1\}$ and
$r_2=\min\{|a_2|: a_2\in A_2\}.$
Then $r_1<r_2$ and the round annulus $R=\{z: r_1<|z|<r_2\}$ separates $A_1$ from $A_2$ and
satisfies the inequality
$$
\mod R=\log\frac{r_2}{r_1}\ge P(A_1,A_2).
$$

Next we show the latter part.
Suppose that a circular annulus $R$ separates $A_1$ from $A_2.$
By a standard limiting process, we may assume, without loss of generality,
that $(A_1\cup A_2)\cap\partial R=\emptyset.$
Fix two pairs of distinct points $a_1,a_1'\in A_1$ and $a_2,a_2'\in A_2\,.$
By using M\"obius invariance of the quantities involved,
we may further assume that $a_2'=\infty$ and that $R$ is of the form $\D\setminus\bD(a,\rho)$
with $0\le a<1-\rho.$
If we set $m=\mod R,$ then $R$ is M\"obius equivalent to the round annulus $r:=e^{-m}<|z|<1.$
Thus we can choose $0\le b<1$ so that $T_b(-r)=a-\rho$ and $T_b(r)=a+\rho,$
where $T_b(z)=(z+b)/(1+bz).$
Since $a=(T_b(r)+T_b(-r))/2$ and $\rho=(T_b(r)-T_b(-r))/2,$ we obtain
$$
\frac{1-a}\rho=\frac{1+br^2}{(1+b)r}>\frac{1+r^2}{2r}.
$$
Noting the inequalities $|a_1-a_2|\ge 1-(a+\rho)$ and $|a_1-a_1'|\le 2\rho,$
we obtain
$$
|\CR(a_1,a_1',a_2,\infty)|=\frac{|a_1-a_2|}{|a_1-a_1'|}
\ge \frac{1-a-\rho}{2\rho}
>\frac{1+r^2}{4r}-\frac12=\frac{(1-r)^2}{4r}=\sinh^2\frac m2.
$$
Hence, $P(A_1,A_2)\ge \log\sinh^2(m/2)=2\log\sinh(m/2).$
We note that $\sinh(m/2)\ge1$ if and only if $m\ge 2\log(1+\sqrt2).$
Since the function $f(x)=2\log\sinh(x/2)-x$ satisfies
$f'(x)=\coth(x/2)-1>0,$ we have $f(x)\ge f(2\log(1+\sqrt2))=-2\log(1+\sqrt2)$
for $x\ge2\log(1+\sqrt2).$
Thus the last inequality has been shown.
\end{pf}

The following result roughly means that $1+Q(A)$ is comparable with $1/\sys(X)$
for $X=\sphere\setminus A.$

\begin{thm}\label{thm:hyp}   
For an arbitrary hyperbolic punctured sphere $X=\sphere\setminus A,$
$$
Q(A)\le\frac{\pi^2}{\sys(X)}.
$$
Conversely, for every $\varepsilon\in(0,\pi^2),$ there exists a constant
$\ell_0=\ell_0(\varepsilon)\in(0,1)$ such that the inequality
$$
\frac{\pi^2-\varepsilon}{\sys(X)}<Q(A)
$$
holds for any hyperbolic punctured sphere $X=\sphere\setminus A$
with $\sys(X)\le\ell_0.$
\end{thm}

\begin{pf}
First we show the first inequality.
We may assume that $Q(A)>0.$
By definition, $Q(A)=P(A_1, A_2)$ for an admissible partition $A=A_1\sqcup A_2.$
By the last lemma, there is a circular annulus $R$ separating $A_1$ from $A_2$
with $\mod R=\log\frac{r_2}{r_1}=P(A_1,A_2)=Q(A).$
We may assume that $R$ is of the form $\{z: r_1<|z|<r_2\}.$
Since the core curve $\gamma: |z|=\sqrt{r_1r_2}$ of $R$ has
the hyperbolic length $\ell_R(\gamma)=\pi^2/\mod R$ (see, for instance,
\cite[(4.2)]{SugawaUP}), we have
$$
Q(A)=\mod R=
\frac{\pi^2}{\ell_R(\gamma)}
\le\frac{\pi^2}{\ell_X(\gamma)}
\le\frac{\pi^2}{\sys(X)},
$$
as required.

We next show the second inequality.
For a given $\varepsilon\in(0,\pi^2),$
choose an $0<\ell_0\le1$ so small that
\begin{equation}\label{eq:lb}
\arctan(1/\sinh\ell_0)>\frac{\pi+\ell_0}2-\frac\varepsilon{2\pi}.
\end{equation}
Assume that $\sys(X)<\ell_0.$
Let $\gamma$ be a closed geodesic of length $\ell<\ell_0$ in $X.$
It suffices to show that $(\pi^2-\varepsilon)/\ell\le Q(A).$
Since $\ell_0\le 1<\arccosh 3=2\log(1+\sqrt2),$
$\gamma$ is a power of a simple closed geodesic
(see \cite[Lemma 7 in p.~235]{Yamada82}).
We may assume that $\gamma$ is simple, namely, $\gamma$ is a Jordan curve
in $X\subset\C.$
Thus $\gamma$ divides $A$ into two parts, say,
$A_1$ and $A_2.$
Since $\gamma$ is homotopic to neither a point in $X$ nor a puncture,
this partition is admissible.
We now show that
$$
P(A_1,A_2)>\frac{\pi^2-\varepsilon}\ell,
$$
which implies the required inequality because $Q(A)\ge P(A_1,A_2).$
Choose two pairs of points $a_1, a_1'\in A_1$ and $a_2, a_2'\in A_2$
so that $P(A_1,A_2)=\log^+|\CR(a_1,a_1',a_2',a_2)|.$
By M\"obius invariance of the cross ratio, we may assume that
$a_1'=0$ and $a_2'=\infty.$
The collar lemma implies that there is an annulus $R$ in $X$
containing $\gamma$ as its core curve with $\mod R>2\pi\arctan(1/\sinh\ell)/\ell$
(Halpern \cite{Hal81}, see also \cite[Theorem 5.2]{SugawaUP}).
By \eqref{eq:lb},  we note that
$$
\mod R>\frac{2\pi\arctan(1/\sinh\ell)}\ell
>\frac{\pi^2-\varepsilon+\pi\ell_0}\ell
>\frac{\pi^2-\varepsilon}\ell+\pi.
$$
A theorem of Avkhadiev-Wirths \cite[Theorem 3.17]{AW:SP} (see also
\cite[\S 3]{Sugawa10JA}) guarantees us the existence of a round annulus $R_0$ of the form
$\{z\in\C: r_1<|z|<r_2\}$ contained in $R$ with
$\mod R_0=\log(r_2/r_1)\ge\mod R-\pi.$
Since $R_0$ separates $A_1$ from $A_2,$ we have $|a_1|\le r_1$ and
$|a_2|\ge r_2,$ and therefore, $|\CR(a_1,0,\infty,a_2)|=\log|a_2/a_1|
\ge\log(r_2/r_1).$
Hence,
$$
\log|\CR(a_1,a_1',a_2',a_2)|\ge \log\frac{r_2}{r_1}\ge\mod R-\pi
>\frac{\pi^2-\varepsilon}\ell
$$
and therefore $P(A_1,A_2)>(\pi^2-\varepsilon)/\ell$ as required.
\end{pf}

We remark that the thrice-punctured sphere $\C_{0,1}=\sphere\setminus\{0,1,\infty\}$
is known to satisfy $\sys(\C_{0,1})=2\log(1+\sqrt2)=1.7627\dots>\ell_0$
(see \cite{Yamada82}).
From the above proof, we obtain also the following.

\begin{cor}
For a punctured spheres $X=\sphere\setminus A$ with $\card(A)\ge 3,$
the inequalities
$$
\frac{C_1}{Q(A)+1}\le\sys(X)\le\frac{C_2}{Q(A)+1}
$$
hold, where $C_1$ and $C_2$ are absolute constants with $1<C_1<C_2< 13.4.$
\end{cor}

\begin{pf}
By the first inequality in Theorem \ref{thm:hyp} and the Schmutz inequality \eqref{eq:Sch},
we obtain
$$
\sys(X)(Q(A)+1)<\pi^2+4\log(1+\sqrt 2)\approx 13.3951.
$$
Thus the right-hand side inequality has been shown.
To show the other one, we take $\ell_0=1$ so that
$$\varepsilon=2\pi[(\pi+\ell_0)/2-\arctan(1/\sinh \ell_0)
=\pi(\pi+1)-2\pi\arctan(1/(e-e\inv))\approx 8.58.$$
If $\sys(X)\le1,$ the last theorem yields
$$
\sys(X)(Q(A)+1)>\sys(X)Q(A)>\pi^2-\epsilon>1.28.
$$
Otherwise, $\sys(X)(Q(A)+1)\ge\sys(X)>1.$
Therefore, we have $\sys(X)(Q(A)+1))>1$ at any event.
\end{pf}

In the inequalities of the last corollary, we cannot replace the denominator
$Q(A)+1$ by $Q(A)$ in general.
Indeed, $Q(A)$ may be very small even when $\sys(X)$ is not so small
as the following example shows.

\begin{example}
Let $X=\sphere\setminus A,$ where $A=\{0,1,a,\infty\}$
for $a\in\C\setminus\{0,1,\infty\}.$
In view of the formula $\CR(0,1,a,\infty)=a,$ we have
$$
Q(A)=\mQ(A)=\max\{|\log|a||,|\log|1-a||,|\log|1-1/a||\}.
$$
We observe that $Q(A)=\mQ(A)\to+\infty$ precisely when
$a$ approaches the boundary of $\C_{0,1}.$
For instance, $Q(A)=\log|a|$ when $a$ lies in
the set $\{z\in\sphere: \Re z\ge1/2, |z-1|\ge 1\}.$
In particular, $Q(A)=0$ when $a=e^{\pm \pi i/3}.$
\end{example}

We next consider a finite subset $A$ of a circle $C$ in $\sphere.$
Suppose that $P(A_1,A_2)>0$ for an admissible partition $A=A_1\sqcup A_2.$
By the former part of Lemma \ref{lem:mod}, we observe that $A_1$
does not separate $A_2$ in $C.$
Therefore, we may assume that $C=\R\cup\{\infty\}$ and
$a_1<a_2$ whenever $a_1\in A_1$ and $a_2\in A_2.$
By \eqref{eq:cr2}, we have
$$
\CR(a_1,a_1',a_2,a_2')=1-\CR(a_1',a_1,a_2,a_2')
=1+\frac{(a_2-a_1')(a_2'-a_1)}{(a_1'-a_1)(a_2'-a_2)}>1
$$
for $a_1,a_1'\in A_1$ and $a_2,a_2'\in A_2$ with $a_1<a_1'$ and $a_2<a_2'.$
Moreover, $\CR(a_1,a_1',a_2,a_2')$ is minimised when $a_2=\min A_2$
and $a_2'=\max A_2$ for a fixed pair $(a_1,a_1').$
Indeed, this can easily be checked by normalising $a_1=0$ and $a_1'=1.$
Similarly, we observe that $\CR(a_1,a_1',a_2,a_2')$ is minimised
when $a_1=\min A_1$ and $a_1'=\max A_1$ for a fixed pair $(a_2,a_2').$
In view of \eqref{eq:cr}, we summarise these observations as follows.

\begin{prop}
Let $A_1$ and $A_2$ be finite subsets of $\R$ such that
$\card(A_j)\ge2,~j=1,2,$ and $a_1<a_2$ for $a_1\in A_1$ and $a_2\in A_2.$
Then
$$
P(A_1,A_2)=\log^+\big\{\CR(\min A_1, \max A_1, \min A_2, \max A_2)-1\big\}.
$$
\end{prop}

The following example shows that $Q(A)$ and $\mQ(A)$ are not comparable.

\begin{example}\label{ex:2}
Let $A=\{0,1,\dots, n\}$ for $n\ge 3.$
If $A_1=\{k+1,k+2,\dots,l-1,l\}$ for some $1\le k+1<l<n,$ then
$\CR(k,l+1,k+1,l)=-1/\{(l-k)^2-1\}$ and thus $P(A_1,A_2)=0.$
If $A_1=\{0,1,\dots,k\}$ for some $1\le k\le n-2,$ then by the above proposition,
$P(A_1,A_2)=\log^+(\CR(0,k,k+1,n)-1)=\log^+\frac{n}{k(n-k-1)},$
which is maximised when $k=1$ or $k=n-2.$
Thus we have $Q(A)=\log\frac n{n-2}.$
On the other hand, since $|(a_3-a_1)(a_4-a_2)|\ge \max\{n(n-2), (n-1)^2\}=(n-1)^2$
and $|(a_2-a_1)(a_4-a_3)|\ge 1$ for any distinct $a_1,a_2,a_3,a_4$ in $A,$
we obtain $|\CR(a_1,a_2,a_3,a_4)|\le (n-1)^2,$ where equality holds
for $(a_1,a_2,a_3,a_4)=(0,1,n-1,n).$
Therefore, $\mQ(A)=2\log(n-1).$
This example shows that $Q(A)\to0$ whereas $\mQ(A)\to+\infty$
as $n\to\infty.$
By Theorem \ref{thm:hyp}, $1+Q(A)$ is comparable with $1/\sys(X)$
but we cannot replace $Q(A)$ by $\mQ(A)$ in the theorem.
\end{example}

\section{Construction of a distance}

Suppose that $f:[0,\pi]\to (0, \infty)$ is a continuous increasing subadditive function  with
$f(0)=0.$
Here, subadditivity means that
$f(t_1+t_2)\le f(t_1)+f(t_2)$ whenever $t_1,t_2,t_1+t_2\in[0,\pi].$
We then define a function $D=D^f$ on $\bDs(0,1/e)^2$ by
\begin{align*}
D^f(z_1,z_2)&=\frac{f(\theta)}{\max\{\log(1/|z_1|),\log(1/|z_2|)\}}
+\left|\log\frac{\log(1/|z_1|)}{\log(1/|z_2|)}\right| \\
&=\frac{f(\theta)}{\max\{\tau_1, \tau_2\}}+|\log\tau_1-\log\tau_2|
\end{align*}
for $z_1,z_2\in\bDs(0,1/e),$
where $\theta=|\arg(z_2/z_1)|=2\arcsin(|\zeta_1-\zeta_2|/2)\in[0,\pi], \zeta_l=z_l/|z_l|,
\tau_l=-\log|z_l|~(l=1,2).$
Note that $D^f(z_1,z_2)=f(\theta)$ for $|z_1|=|z_2|=e^{-1}.$
Our construction will be based on the following lemma.

\begin{lem}\label{lem:Df}
Let $L_1$ and $L_2$ be positive numbers and let $f:[0,\pi]\to (0, \infty)$ be a continuous increasing subadditive function
satisfying $L_1t\le f(t)\le L_2t$ for $0\le t\le \pi.$
Then
$$
2L_1 h_{\D^*}(z_1,z_2)\le D^f(z_1,z_2)\le M_0L_2 h_{\C_{0,1}}(z_1,z_2),
\quad 0<|z_1|,|z_2|\le e^{-1},
$$
where $M_0$ is a constant with $M_0<24.$
If $f\le 2$ in addition, then $D^f(z_1,z_2)$ satisfies the triangle inequality
on $0<|z|\le e\inv.$
\end{lem}

The first part is proved in \cite{SZ17} for $f(t)=2\sin(t/2)$ and $f(t)=\id(t)=t.$
The second part is also proved for $f(t)=2\sin(t/2)$ in \cite{SZ17}.
We include an outline of the proof here for the convenience of the reader.
For details, see \cite{SZ17}.

For instance, if $f$ is
a continuously differentiable, increasing, and concave function in $[0,\pi]$
with $f(0)=0,$ then $f(t)/t$ is non-increasing and $f$ is subadditive.
In particular, the inequalities
$$
\frac{f(\pi)}{\pi}\le \frac{f(t)}{t}\le f'(0),\quad 0<t\le \pi,
$$
hold. Thus we can take $L_1=f(\pi)/\pi$ and $L_2=f'(0)$ in this case.

\begin{pf}[Proof of Lemma \ref{lem:Df}]
For the first part, we note that the inequalities
$L_1 D^\id(z_1,z_2)\le D^f(z_1,z_2)\le L_2 D^\id(z_1,z_2)$ hold.
Since the inequalities $2h_{\D^*}(z_1,z_2)\le D^\id(z_1,z_2)\le M_0h_{\C_{0,1}}(z_1,z_2)$
are already shown in \cite{SZ17}, the required inequalities follow.

We next show the triangle inequality
$\Delta:=D(z_1,z)+D(z,z_2)-D(z_1,z_2)\ge0$
under the assumption $f\le2.$
The most delicate case is when $1\le\tau_2\le\tau_1\le\tau,$
where $\tau_l=-\log|z_l|~(l=1,2)$ and $\tau=-\log|z|.$
The other cases can be handled similarly or even more easily, and thus,
will be omitted.
Let $\theta=|\arg(z_2/z_1)|, \theta_l=|\arg(z_l/z)|\in[0,\pi]~(l=1,2).$
Since $\theta\le \theta_1+\theta_2,$ it is easy to check that the inequality
$f(\theta)\le f(\theta_1)+f(\theta_2)$ holds.
Thus,
\begin{align*}
\Delta
&=\frac{f(\theta_1)+f(\theta_2)}{\tau}+\log\frac{\tau}{\tau_1}+\log\frac\tau{\tau_2}
-\frac{f(\theta)}{\tau_1}-\log\frac{\tau_1}{\tau_2} \\
&\ge\frac{f(\theta)}{\tau}-\frac{f(\theta)}{\tau_1}+2\log\tau-2\log\tau_1
=: g(\tau).
\end{align*}
Since $g'(\tau)=[2\tau-f(\theta)]/\tau^2\ge[2-f(\theta)]/\tau^2\ge0,$
one has $\Delta\ge g(\tau)\ge g(\tau_1)=0.$

\end{pf}

\begin{rem}
We remark that if ``maximum" is replaced by ``minimum" in the definition of
the function $D^f(z_1,z_2),$ then the corresponding quantity does not necessarily
satisfy the triangle inequality.
\end{rem}

We enumerate the elements of $A$ as $a_1,\dots, a_n$ so that
$
a_n=\infty.$
Set
$$
\tilde\rho_j=\begin{cases}
\displaystyle
\min_{1\le k<n, k\ne j}|a_k-a_j|&\quad \text{for}~j=1,2,\dots, n-1, \\
\displaystyle
\max_{1\le k<n}|a_k| &\quad \text{for}~j=n,
\end{cases}
$$
$\rho_j=\tilde\rho_j/e$ for $1\le j<n$ and $\rho_n=e\tilde\rho_n.$
We next consider the sets
$E_j=\bD(a_j,\rho_j)$ for $1\le j<n$ and $E_n=\{z\in\sphere:|z|\ge\rho_n\}.$
Furthermore, set
$E_j^*=E_j\setminus\{a_j\}$
and $W=\sphere\setminus(E_1\cup\dots\cup E_n).$
Note that $E_1, \dots, E_n$ are disjoint.
Since $0\in A,$ we observe that $e\rho_j=\tilde\rho_j\le |a_j|\le
\tilde\rho_n=\rho_n/e$ for $1\le j<n.$
In particular,
\begin{equation}\label{eq:max}
\max_{1\le j<n}\rho_j\le e^{-2}\rho_n.
\end{equation}
We set
$$
\rho_{\min}=\min_{1\le j\le n}\rho_j
\aand
\rho_{\max}=\max_{1\le j\le n}\rho_j=\rho_n.
$$

Suppose that $\mu(z_1,z_2)$ is a distance function on $X$ such that
\begin{equation}\label{eq:mu}
\mu(z_1,z_2)=S_jf_{j}(\theta) \quad \text{for}~ z_1,z_2\in\partial E_j, ~j=1,\dots,n,
\end{equation}
where $\theta=|\arg((z_2-a_j)/(z_1-a_j))|\in[0,\pi]$ for $j<n$ and
$\theta=|\arg (z_2/z_1)|\in[0,\pi]$ for $j=n,$
$f_j:[0,\pi]\to[0,2]$ is a continuous increasing subadditive function
with $f_j(0)=0$ and $S_j$ is a positive constant for each $j.$
We note that by \eqref{eq:mu} the function $\mu(z_1,z_2)$ is rotationally invariant on $\partial E_j$
about $a_j.$

\begin{example}\label{ex:Euc}
When $\mu(z_1,z_2)$ is the euclidean distance $|z_1-z_2|,$
we can choose $f_j(t)=2\sin(t/2)$ and $S_j=\rho_j$ in the above construction.
\end{example}

Then, we define a function $D_j(z_1,z_2)$ on $E_j^*\times E_j^*$ by
$$
D_{j}(z_1,z_2)=
\begin{cases}
\displaystyle
S_j D^{f_j}\left(\frac{z_1-a_j}{\tilde\rho_j},\frac{z_2-a_j}{\tilde\rho_j}\right)
&\quad\text{if}~1\le j<n, \\
\null & \null \\
\displaystyle
S_n D^{f_n}\left(\frac{\tilde\rho_n}{z_1},\frac{\tilde\rho_n}{z_2}\right)
&\quad\text{if}~ j=n,
\end{cases}
$$
for $z_1, z_2\in E_j^*.$
Then $D_j(z_1,z_2)=\mu(z_1,z_2)$ for $z_1,z_2\in\partial E_j.$
We now apply the following lemma with the choice $Y_0=W,~ d_0=\mu,~
Y_j=E_j^*,~d_j=D_j~(1\le j\le n)$ to construct a new distance $d$ on $X.$
We will denote by $d_X$ this distance $d$ constructed in this way with $\mu, f_j, S_j~(j=1,\dots,n).$
The same construction was given in \cite{SZ17} when $\mu(z_1,z_2)=|z_1-z_2|.$

\begin{lem}\label{lem:dist}
Let $X$ be a set and assume that $X=Y_0\cup Y_1\cup\cdots\cup Y_n.$
Let $d_j$ be a distance function on $Y_j$ for $j=0,1,\dots,n.$
Suppose that $Y_j\cap Y_k=\emptyset$ for $j\ne k$ with $1\le j,k\le n,$
that $d_0(x_1,x_2)=d_j(x_1,x_2)$ for $x_1,x_2\in Z_j=Y_0\cap Y_j,~j=1,\dots, n,$ and that
each $Z_j$ is non-empty and closed in the both metric spaces $(Y_0,d_0)$ and $(Y_j,d_j).$
We define $d(x_1,x_2)$ for $x_1,x_2\in X$ to be
$$
\begin{cases}
d_j(x_1,x_2)
&\quad\text{for}~ x_1,x_2\in Y_j' ~\text{with}~0\le j\le n, \\
\displaystyle
\inf_{\zeta\in Z_j}\big\{ d_j(x_1,\zeta)+d_0(\zeta,x_2)\big\}
&\quad\text{for}~ x_1\in Y_j, x_2\in Y_0~\text{with}~j\ne0, \\
\displaystyle
\inf_{\zeta\in Z_j}\big\{d_0(x_1,\zeta)+d_j(\zeta,x_2)\big\}
&\quad\text{for}~ x_1\in Y_0, x_2\in Y_j~\text{with}~j\ne0, \\
\displaystyle
\inf_{\substack{\zeta_1\in Z_j \\ \zeta_2\in Z_k}}
\big\{ d_j(x_1,\zeta_1)+d_0(\zeta_1,\zeta_2)+d_k(\zeta_2,x_2)\big\}
&\quad\text{for}~ x_1\in Y_j, x_2\in Y_k~\text{with}~ 0\ne j\ne k\ne 0,
\end{cases}
$$
where $Y_j'=Y_j\setminus Z_j$ for $j=1,\dots,n$ and $Y_0'=Y_0\setminus(Z_1\cup\dots\cup Z_n).$
Then $d$ is a distance function on $X$ such that $d(x_1,x_2)=d_j(x_1,x_2)$
for $x_1,x_2\in Y_j,~j=0,1,\dots,n.$
\end{lem}

\begin{pf}
First we check that $d=d_j$ on $Y_j\times Y_j.$
For instance, we assume that $x_1\in Y_j'$ and $x_2\in Z_j$ for some $j=1,\dots,n.$
Letting $\zeta=x_2$ in the second case of the definition of $d(x_1,x_2),$ we see that
$d(x_1,x_2)\le d_j(x_1,x_2).$
On the other hand, for any $\zeta\in Z_j,$ by the triangle inequality and the assumption, we have
$$
d_j(x_1,\zeta)+d_0(\zeta,x_2)
=d_j(x_1,\zeta)+d_j(\zeta,x_2)\ge d_j(x_1,x_2),
$$
which implies $d(x_1,x_2)\ge d_j(x_1,x_2).$
Thus we have $d(x_1,x_2)=d_j(x_1,x_2).$
The other cases can be handled similarly.
We thus checked the last assertion.

We next show that $d$ is a distance function on $X.$
The symmetricity $d(x_1,x_2)=d(x_2,x_1)$ and nonnegativity are obvious by definition.
We now show that $d(x_1,x_2)=0$ only if $x_1=x_2.$
Suppose that $d(x_1,x_2)=0$ for a pair of points $x_1,x_2$ in $X.$
When $x_1$ and $x_2$ are contained in a common $Y_j,$ then $d(x_1,x_2)=d_j(x_1,x_2)=0$
implies $x_1=x_2$ certainly.
We now consider, for example, the case when $x_1\in Y_j'$ and $x_2\in Y_0$ for some $j\ne0.$
Then, by definition,
$$
0=d(x_1,x_2)\ge \inf_{\zeta\in Z_j}d_j(x_1,\zeta).
$$
This implies that $x_1$ is contained in the closure of $Z_j,$
which is a contradiction because $Z_j$ is closed in $Y_j.$
The other cases can be excluded in a similar way.
Hence, we conclude that $d(x_1,x_2)=0$ implies $x_1=x_2.$

We finally check the triangle inequality.
As a sample case, we assume that $x_1\in Y_j, x_2\in Y_0, x_3\in Y_k$
with $j\ne0\ne k\ne j.$
Then
\begin{align*}
d(x_1,x_2)+d(x_2,x_3)
&= \inf_{\zeta_1\in Z_j,\zeta_2\in Z_k}\big\{ d_j(x_1,\zeta_1)+d_0(\zeta_1,x_2)+d_0(x_2,\zeta_2)+d_k(\zeta_2,x_3)\big\} \\
&\ge \inf_{\zeta_1\in Z_j,\zeta_2\in Z_k}\big\{d_j(x_1,\zeta_1)+d_0(\zeta_1,\zeta_2)+d_k(\zeta_2,x_3)\big\} \\
&=d(x_1,x_3).
\end{align*}
We can deal with the other cases similarly.
\end{pf}

We now state a general result.
This was already shown in \cite{SZ17} under the situation of Example \ref{ex:Euc}.
Indeed, when $\mu(z_1,z_2)=|z_1-z_2|,$
a similar quantity $e_X(z_1,z_2)$ is defined without taking the infima
and a comparison result is given for it as Theorem 5.1 in  \cite{SZ17}.
Since its proof is tedious,
we restrict our attention to $d_X(z_1,z_2)$ for simplicity in the present note.

\begin{thm}\label{thm:general}
Let $X=\sphere\setminus A$ be an $n$-times punctured sphere with $n\ge 3$
and $X\subset\C\setminus\{0,1\}.$
Let $d_X$ be the distance on $X$ constructed in Lemma \ref{lem:dist} with a distance $\mu(z_1,z_2)$
satisfying the condition \eqref{eq:mu} with $f_j:[0,\pi]\to[0,2], S_j>0, ~(j=1,\dots,n).$
Suppose that $L_1 t\le f_j(t)\le L_2t$ for $0\le t\le\pi$ and $j=1,\dots,n$
and that $K_1 h_X(z_1,z_2)\le \mu(z_1,z_2)\le K_2h_X(z_1,z_2)$
for $z_1,z_2\in W=X\setminus(E^*_1\cup\cdots\cup E^*_n).$
Then the distance $d_X(z_1,z_2)$ on $X$ satisfies the inequalities
$$
B_1 h_X(z_1,z_2)\le d_X(z_1,z_2) \le B_2 h_X(z_1,z_2),\quad z_1,z_2\in X,
$$
where
$$
B_1=\min\{2S_{\min}L_1,K_1\}
\aand
B_2=\max\{S_{\max} M_0L_2, K_2\}
$$
and $S_{\max}=\max\{S_1,\dots,S_n\},~S_{\min}=\min\{S_1,\dots,S_n\},~ M_0<24.$
\end{thm}

\begin{pf}
Since the method is almost the same as in \cite[Theorem 5.1]{SZ17}, we give only a sketchy proof.
As a typical case, we consider the situation that $z_1\in E_j^*$ and $z_2\in W.$
Let $\alpha$ be a shortest hyperbolic geodesic joining $z_1$ and $z_2$ in $X.$
Then $h_X(z_1,z_2)=h_X(z_1,\zeta_0)+h_X(\zeta_0,z_2)$ for the first hitting point
$\zeta_0$ of $\alpha$ to $\partial E_j$ starting from $z_1.$
We first assume that $j\ne n.$
Choose $a_k$ so that $\tilde\rho_j=|a_k-a_j|$
and set $X_1=\C\setminus\{a_j,a_k\}.$
By Lemma \ref{lem:Df}, we obtain
$$
D_j(z_1,\zeta_0)\le S_jM_0L_2 h_{X_1}(z_1,\zeta_0) \le S_jM_0L_2 h_{X}(z_1,\zeta_0).
$$
Next we assume that $j=n.$
Choose $a_k$ so that $\tilde\rho_n=|a_k|$ and set $g(z)=a_k/z.$
Then, by Lemma \ref{lem:Df} again, we have
$$
D_n(z_1,\zeta_0)\le S_nM_0L_2h_{\C_{0,1}}(g(z_1),g(\zeta_0))\le S_nM_0L_2h_X(z_1,\zeta_0).
$$
Since $\mu(\zeta_0,z_2)\le K_2h_X(\zeta_0,z_2)$ by assumption, we obtain
$$
d_X(z_1,z_2)
\le D_j(z_1,\zeta_0)+\mu(\zeta_0,z_2)
\le B_2h_X(z_1,z_2).
$$

On the other hand, by definition,
$d_X(z_1,z_2)=D_j(z_1,\zeta)+\mu(\zeta,z_2)$ for some $\zeta\in\partial E_j.$
By Lemma \ref{lem:Df}, in a similar way as above, we obtain
$D_j(z_1,\zeta)\ge 2S_jL_1h_X(z_1,\zeta).$
Thus,
$$
d_X(z_1,z_2)\ge 2S_{\min} L_1h_X(z_1,\zeta)+K_1h_X(\zeta,z_2)\ge B_1 h_X(z_1,z_2).
$$
Thus the proof is complete.
\end{pf}

As an example, we apply the last theorem to Example \ref{ex:Euc}.
We note here that one can choose $L_1=2/\pi,~L_2=1, K_1=2\rho_{\min}/\pi$ and
\begin{equation}\label{eq:K}
K_2=\max\{2\rho_ne^{-K_0}(C_0-K_0-1+\log[\rho_n/\rho_{\min}]),
6\rho_n(C_0+\log[3\rho_n])\},
\end{equation}
where $C_0=1/(2\lambda_{\C_{0,1}}(-1))=\Gamma(1/4)^4/4\pi^2\approx 4.37688$
and $K_0\approx 0.846666$ (see the proof of Theorem 5.1 in \cite{SZ17} for details).
Thus, in view of \eqref{eq:max},
we obtain the following result, which is essentially contained in \cite{SZ17}.

\begin{cor}\label{cor:Euc}
Let $X=\sphere\setminus \{a_1,\dots,a_n\}$ be an $n$-times punctured sphere with
$a_1=0, ~ a_n=\infty$ and  let $d_X$ be the distance
on $X$ constructed above for the choices
$\mu(z_1,z_2)=|z_1-z_2|, f_j(t)=2\sin(t/2)$ and $S_j=\rho_j$ for $j=1,\dots,n.$
Then
$$
B_1 h_X(z_1,z_2)\le d_X(z_1,z_2) \le B_2 h_X(z_1,z_2),\quad z_1,z_2\in X,
$$
where
$$
B_1=\frac{2\rho_{\min}}\pi
\aand
B_2=\max\{M_0\rho_n, K_2\},
$$
$K_2$ is given in \eqref{eq:K}
and $\rho_{\min}=\min\{\rho_1,\dots,\rho_n\}\,.$
\end{cor}

\section{quasihyperbolic type constructions}
Let $X=\sphere\setminus A$ be an $n$-times punctured sphere as before,
where $A$ is a finite subset of $\sphere$ with $0,1,\infty\in A.$
We recall that $\delta_X(z)$ denotes the euclidean distance from $z\in X$
to the boundary $\partial X=A.$ For all   $z$ in an open dense subset of $ X\,,$ there is only one $a\in A$ such that
$\delta_X(z)=|z-a|.$
For $a\in A,$ we set $V_a=\{z\in\C: |z-a|\le|z-b| ~\text{for any}~ b\in A \}.$
Note that each $V_a$ is a closed convex set and that $V_\infty=\emptyset.$
Then we obtain a partition of $\C$ into $V_a,~ a\in A,$ which is known as
the {\it Voronoi diagram}.
Each $V_a$ is called a Voronoi cell and the point $a$ is called the nucleus of $V_a.$
Since $\C$ is tessellated by $V_a$'s, some of them are unbounded.
It is known that $V_a$ is unbounded if and only if $a$ lies on the boundary of
the convex hull of $A\setminus\{\infty\}.$
See Figure 1 (a).
We refer to \cite{Vai09} for the relation between the quasihyperbolic metric and the Voronoi diagram.
Set $\delta_a(z)=|z-a|.$
Then $\delta_X(z)=\min_{a\in A}\delta_a(z).$
Since $\delta_X(z)=\delta_a(z)$ for $z\in V_a,$ we see that the function $\delta_X(re^{i\theta})$
is not of class $C^1$ in $\theta\in\R$ (and thus, not rotationally invariant) for large enough $r>0.$
Therefore, the $j$-distance (see \cite{AVV:conf} for instance)
$$
j_X(z_1,z_2)=\log\left(1+\frac{|z_1-z_2|}{\min\{\delta_X(z_1),\delta_X(z_2)\}}\right)
$$
does not satisfy the condition \eqref{eq:mu}.

Therefore, we are led to a slight modification of the function $\delta_X$ in the following way:
$$
\deltat_X(z)=\min\{\delta_X(z), |z|/2\}.
$$
Here we recall that the set $A$ of the punctures contains $0$ and $\infty.$
Note that $\delta_X(z)/2\le \deltat_X(z)\le\delta_X(z)$ for $z\in X.$
In other words, $\deltat_X(z)=\min_{a\in A}\deltat_a(z),$ where
$$
\deltat_a(z)=\begin{cases}
|z-a| & \text{if}~ a\ne 0, \\
|z|/2 & \text{if}~ a=0, \\
+\infty & \text{if}~a=\infty.
\end{cases}
$$
Then the {\it modified Voronoi diagram} (with respect to the origin) is defined by
setting
$$
\Vt_a=\{z\in\C: \deltat_a(z)\le \deltat_b(z) ~\text{for any}~ b\in A\}
$$
as the modified Voronoi cell for $a\in A.$
For $a,b\ne0,$ the set $H_{a,b}=\{z: \deltat_a(z)\le \deltat_b(z)\}$ is the half-plane
containing $a$  and bounded by the bisector of the segment joining the two points $a$ and $b.$
It is easy to see that $H_{a,0}$ is the Apollonian disk $\{z: |z-a|\le |z|/2\};$ namely,
$\{z: |z-4a/3|\le 2|a|/3\}.$
Therefore, for $a\in A\setminus\{0,\infty\},$
$\Vt_a$ is the intersection of a (possibly unbounded) convex polygon with
the disk $H_{a,0}$ and, in particular, a compact convex set.
On the other hand, $\Vt_0$ is a (possibly disconnected, non-convex) unbounded closed set.
By definition, we have $\deltat_X(z)=\deltat_a(z)$ for $z\in\Vt_a,~a\in A.$
See Figure 1 (b) for the modified Voronoi diagram.

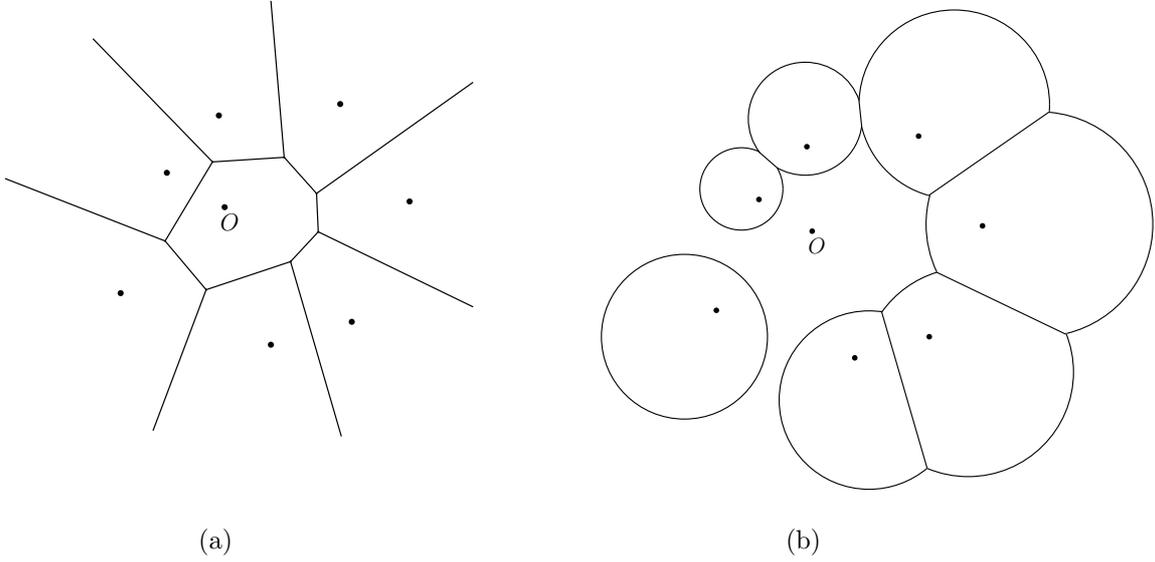
\begin{figure}[htbp]
\centering
 \hspace{-19mm}\subfigure[\label{fig1}]{\scalebox{0.76}[0.76]{
\begin{pspicture}(-3.3,-5.3)(3.3,3.3)

\psdot[dotsize=3pt](0,0)
\psdot[dotsize=3pt](-0.1,1.6)
\psdot[dotsize=3pt](3.2,0.1)
\psdot[dotsize=3pt](2,1.8)
\psdot[dotsize=3pt](2.2,-2)
\psdot[dotsize=3pt](0.8,-2.4)
\psdot[dotsize=3pt](-1.8,-1.5)
\psdot[dotsize=3pt](-1,0.6)

\psdot[dotsize=1pt](-0.21,0.79)
\psdot[dotsize=1pt](1.03,0.87)
\psdot[dotsize=1pt](1.59,0.24)
\psdot[dotsize=1pt](1.62,-0.43)
\psdot[dotsize=1pt](1.14,-0.95)
\psdot[dotsize=1pt](-0.32,-1.44)
\psdot[dotsize=1pt](-1.03,-0.59)

\psline[linewidth=0.6pt](-0.21,0.79)(1.03,0.87)
\psline[linewidth=0.6pt](1.59,0.24)(1.03,0.87)
\psline[linewidth=0.6pt](1.59,0.24)(1.62,-0.43)
\psline[linewidth=0.6pt](1.62,-0.43)(1.14,-0.95)
\psline[linewidth=0.6pt](1.14,-0.95)(-0.32,-1.44)
\psline[linewidth=0.6pt](-0.32,-1.44)(-1.03,-0.59)
\psline[linewidth=0.6pt](-1.03,-0.59)(-0.21,0.79)


\psline[linewidth=0.6pt](-0.21,0.79)(-2.28,2.94)
\psline[linewidth=0.6pt](1.03,0.87)(0.8,3.6)
\psline[linewidth=0.6pt](1.59,0.24)(4.3,2.18)
\psline[linewidth=0.6pt](1.62,-0.43)(4.3,-1.74)
\psline[linewidth=0.6pt](1.14,-0.95)(2.02,-4)
\psline[linewidth=0.6pt](-0.32,-1.44)(-1.24,-3.9)
\psline[linewidth=0.6pt](-1.03,-0.59)(-3.8,0.5)

\uput[d](0.09,0.05){$O$}

\end{pspicture}}}\ \ \hspace{24mm}
  \subfigure[\label{fig2}]{\scalebox{0.7}[0.7]{
\begin{pspicture}(-3.3,-5.3)(3.3,3.3)

\psdot[dotsize=3pt](0,0)
\psdot[dotsize=3pt](-0.1,1.6)
\psdot[dotsize=3pt](3.2,0.1)
\psdot[dotsize=3pt](2,1.8)
\psdot[dotsize=3pt](2.2,-2)
\psdot[dotsize=3pt](0.8,-2.4)
\psdot[dotsize=3pt](-1.8,-1.5)
\psdot[dotsize=3pt](-1,0.6)

\psarc[linewidth=0.6pt](-1.33,0.8){0.78}{0}{32}
\psarc[linewidth=0.6pt](-1.33,0.8){0.78}{65}{360}

\psarc[linewidth=0.6pt](-0.13,2.13){1.07}{19}{217}
\psarc[linewidth=0.6pt](-0.13,2.13){1.07}{240}{352}

\psarc[linewidth=0.6pt](2.67,2.4){1.79}{356}{360}
\psarc[linewidth=0.6pt](2.67,2.4){1.79}{0}{178}
\psarc[linewidth=0.6pt](2.67,2.4){1.79}{193}{254}

\psarc[linewidth=0.6pt](4.27,0.13){2.13}{0}{85}
\psarc[linewidth=0.6pt](4.27,0.13){2.13}{165}{205}
\psarc[linewidth=0.6pt](4.27,0.13){2.13}{284}{360}

\psarc[linewidth=0.6pt](2.93,-2.67){1.98}{0}{22}
\psarc[linewidth=0.6pt](2.93,-2.67){1.98}{107}{145}
\psarc[linewidth=0.6pt](2.93,-2.67){1.98}{247}{360}

\psarc[linewidth=0.6pt](1.07,-3.2){1.69}{82}{310}

\psarc[linewidth=0.6pt](-2.4,-2){1.56}{0}{360}

\psline[linewidth=0.6pt](-0.67,1.21)(-1,1.5)
\psline[linewidth=0.6pt](0.93,1.99)(0.88,2.47)
\psline[linewidth=0.6pt](2.18,0.67)(4.47,2.27)
\psline[linewidth=0.6pt](2.34,-0.78)(4.76,-1.95)
\psline[linewidth=0.6pt](1.3,-1.51)(2.16,-4.49)

\uput[d](0.09,0.05){$O$}

\end{pspicture}}}
\vspace*{-3mm}
\caption{Voronoi diagram (a) and modified Voronoi diagram (b)}
\end{figure}

We enumerate $A$ as $a_1,\dots, a_n$ so that $a_1=0$ and $a_n=\infty,$ as before.
Then we have the next result.

\begin{lem}\label{lem:j}
$\deltat_X(z)=|z|/2$ if either $z\in E_1^*$ or $z\in E_n^*.$
For $z\in E_j^*$ with $1<j<n,$ $\deltat_X(z)=|z-a_j|.$
\end{lem}

\begin{pf}
Let $z\in E_n$ and $a$ be a finite point in $A.$
Since $|a|\le 1/\tilde\rho_n=1/e\rho_n,$ we obtain
$|z-a|\ge|z|-1/e\rho_n\ge |z|/2.$
Thus $\delta_X(z)\ge |z|/2$ for $z\in E_n.$
The other cases can also be checked similarly.
\end{pf}

Applying this idea, we define a modified version of the $j$-distance as follows
$$
\jt_X(z_1,z_2)=\log\left(1+\frac{|z_1-z_2|}{\min\{\deltat_X(z_1),\deltat_X(z_2)\}}\right) \,.
$$
Noting the elementary inequality $|\delta_X(z_1)-\delta_X(z_2)|
\le|z_1-z_2|,$ the following assertion can be verified in
the same way as in the case of the $j$-metric
\cite[Thm 7.47, p.~147]{AVV:conf}, \cite[p.~36]{GH:qd}.

\begin{lem} \label{hatj}
The function $\jt_X : X\times X \to [0,\infty)$ satisfies the triangle inequality. Moreover,
$  1 \le \jt_X(x,y)/j_X(x,y) \le 2$ for all $x,y \in X \,.$
\end{lem}

In view of Lemma \ref{hatj}, $\jt_X$ is a distance on $X$ which is comparable with the $j$-metric $j_X.$
By Lemma \ref{lem:j}, for $z_1,z_2\in \partial E_j$ with $1<j<n,$
$\deltat_X(z)=|z-a_j|=\rho_j$ and thus
$j_X(z_1,z_2)=\log(1+|z_1-z_2|/\rho_j).$
Similarly, we have $j_X(z_1,z_2)=\log(1+2|z_1-z_2|/\rho_j)$
for $z_1,z_2\in\partial E_j$ with $j=1,n.$
If we take $f_j(t)=\log(1+2\sin(t/2))$ for $1<j<n$ and
$f_j(t)=\log(1+4\sin(t/2))$ for $j=1,n$ and $S_j=1$ for $1\le j\le n$ in the above construction,
we obtain the relation $D_j(z_1,z_2)=\jt_X(z_1,z_2)$ for
$z_1,z_2\in \partial E_j,~1\le j\le n.$
Since $\jt(z_1,z_2)$ is Lipschitz equivalent to the euclidean distance $|z_1-z_2|$
on $W,$ the next result follows from Theorem \ref{thm:general}.

\begin{thm}\label{thm:j}
Suppose that $A$ is a subset of $\sphere$ consisting of $n\ge3$ points with
$0,\infty\in A$ and $X=\sphere\setminus A.$
The distance $d_X$ constructed above with $\mu(z_1,z_2)=\jt_X(z_1,z_2)$
satisfies the inequalities
$$
B_1 d_X(z_1,z_2)\le h_X(z_1,z_2) \le B_2 d_X(z_1,z_2),\quad z_1,z_2\in X,
$$
where $B_1$ and $B_2$ are positive constants depending on $X.$
\end{thm}

Similarly, we can modify the definition of  the quasihyperbolic distance by setting
$$
\qt_X(z_1,z_2)=\inf_{\gamma}\int_\gamma\frac{|dw|}{\deltat_X(z)},
$$
where the infimum is taken over all the rectifiable curves $\gamma$ joining
$z_1$ and $z_2$ in $X.$
In the proof of \eqref{eq:MO}, Martin and Osgood indeed showed that
a quasihyperbolic geodesic in $\C^*$ is an arc of a logarithmic spiral about the origin.
In particular, the punctured disk $0<|z|<r$ is quasihyperbolically convex.
Since $1/\delta_X(z)\le 1/\deltat_X(z),$ Lemma \ref{lem:j} implies that
$E_j^*$ is convex with respect to the metric $\qt_X$ and the following
representation is valid:
\begin{equation}\label{eq:qt}
\qt_X(z_1,z_2)=\begin{cases}
q_X(z_1,z_2)=q(z_1-a_j,z_2-a_j), & z_1,z_2\in E_j^*,~ 1<j<n, \\
q(z_1,z_2)/2, & z_1,z_2\in E_j^*,~ j=1,n.
\end{cases}
\end{equation}
Here, $q(z_1,z_2)$ is given in \eqref{eq:MO}.
In particular, we obtain
$$
\qt_X(z_1,z_2)=\begin{cases}
\theta, & z_1,z_2\in\partial E_j^*,~ 1<j<n, \\
\theta/2, & z_1,z_2\in\partial E_j^*,~ j=1,n,
\end{cases}
$$
where $\theta=|\arg((z_2-a_j)/(z_1-a_j))|\in[0,\pi]$ for $1\le j<n$
and $\theta=|\arg (z_2/z_1)|\in[0,\pi]$ for $j=n.$
Therefore, in the construction of $d_X$ in the previous section, we can choose
$f_j(t)=(2/\pi)t$ for $j=1,\dots,n,$ $S_j=\pi/2$ for $1<j<n$
and $S_j=\pi/4$ for $j=1,n.$

\begin{thm}\label{thm:q}
Let $X=\sphere\setminus A$ be an $n$-times punctured sphere with
$X\subset\C_{0,1}.$
Then the distance $d_X(z_1,z_2)$ on $X$ constructed above with
$\mu(z_1,z_2)=\qt_X(z_1,z_2)$ is comparable with the hyperbolic metric on $X:$
$$
B_1 d_X(z_1,z_2)\le h_X(z_1,z_2) \le B_2 d_X(z_1,z_2),\quad z_1,z_2\in X,
$$
where $B_1$ and $B_2$ are positive constants depending only on $Q(A).$
\end{thm}

In Corollary \ref{cor:Euc} and Theorem \ref{thm:j}, the constants $B_1$ and $B_2$
may depend on individual $X$ (more precisely, on $\rho_{\min}$ and $\rho_{\max}=\rho_n$).
In view of Proposition \ref{prop:comp}, we may choose $B_1$ and $B_2$
as a function of $\mQ(A).$
The merit of this theorem is that the constants depend on $Q(A)$ only.

We will need general estimates of the hyperbolic metric.
Following Beardon and Pommerenke \cite{BP78},
we define the technical but useful quantity
$$
\beta_\Omega(z)=\inf\left\{\left|\log\left|\frac{z-a}{b-a}\right|\right|:
a,b\in\partial\Omega, |z-a|=\delta_\Omega(z)\right\}
$$
for $z\in \Omega\subset\C.$
Note that $\beta_\Omega(z)=\beta_{f(\Omega)}(f(z))$ for an affine
automorphism $f(z)=kz+l,~k\ne0, $ of $\C.$
The following result describes nicely the behaviour of the hyperbolic metric.

\begin{lem}[Beardon-Pommerenke \cite{BP78}]\label{lem:BP}
$$
C_2\le\lambda_\Omega(z)\delta_\Omega(z)(C_1+\beta_\Omega(z))\le C_3,
$$
where $C_1, C_2, C_3$ are positive absolute constants.
\end{lem}

Indeed, Beardon and Pommerenke \cite{BP78} gave $C_1=4+\log(3+2\sqrt2), C_2=1/2\sqrt2, C_3=C_1+\pi/4.$
The proof of Theorem \ref{thm:q} will be based on the following lemma.

\begin{lem}\label{lem:beta}
Let $A$ be a finite subset of $\sphere$ with $n=\card(A)\ge3$ and $\infty\in A.$
Then $X=\sphere\setminus A$ satisfies the inequality
$$
\beta_X(z)\le \max\big\{\tfrac12Q(A)+c, \log(1+2e)\big\}, \quad z\in W,
$$
where $W$ is as in  Theorem \ref{thm:general} and $c=\frac12\log[1+1/(4e(1+e))]\approx 0.012217.$
\end{lem}

\begin{pf}
Put $Q=Q(A).$
We may assume that $0\in A$ because of the affine invariance of
the quantities involved.
As before, we enumerate $A$ as $a_1,\dots,a_n$ so that $a_n=\infty$
(but we do not assume $a_1=0$ at this time).
On the contrary, suppose that
$\beta_0:=\beta_X(z_0)>\max\{Q/2+c, \log(1+2e)\}$
for some point $z_0\in W.$
Note here that $\log(1+2e)\approx 1.86199.$
Then $\delta_0:=\delta_X(z_0)=|z_0-a|$
for some $a\in A.$
We assume $a=a_1.$
By the definition of $\beta_X(z_0),$ the set $A$ decomposes into the two sets
$A_1=\{a_j: |a_j-a_1|\le e^{-\beta_0}\delta_0\}$ and
$A_2=\{a_j: |a_j-a_1|\ge e^{\beta_0}\delta_0\}.$
Suppose that $\card(A_1)\ge2$ and $\card(A_2)\ge 2.$
Then $A=A_1\sqcup A_2$ and thus
$$
|\CR(a_1,a_j,a_k,a_n)|=\frac{|a_k-a_1|}{|a_j-a_1|}\ge e^{2\beta_0}
$$
for $a_j\in A_1\setminus\{a_1\}, a_k\in A_2\setminus\{a_n\}.$
In view of \eqref{eq:cr2}, we have thus
$P(A_1,A_2)\ge \log(e^{2\beta_0}-1)=2\beta_0-f(2\beta_0),$
where $f(x)=-\log(1-e^{-x}).$
Since $f(x)$ is decreasing in $0<x<+\infty,$ we observe that
$f(2\beta_0)<f(2\log(1+2e))=\log\frac{1+4e+4e^2}{4e+4e^2}=2c.$
Hence, $Q(A)\ge P(A_1,A_2)>2\beta_0-2c,$ which is a contradiction.

Therefore, we should have either
\begin{enumerate}
\item[(i)]
$|a_j-a_1|\le e^{-\beta_0}\delta_0$ for all $1<j<n,$ or else
\item[(ii)]
$|a_j-a_1|\ge e^{\beta_0}\delta_0$ for all $1<j<n.$
\end{enumerate}
In the case (ii), we have $|z_0-a_1|=\delta_0\le e^{-\beta_0}\tilde\rho_1
<\rho_1$ because $\beta_0>1.$
Therefore, $z_0\in E_1,$ which violates the assumption $z_0\in W.$
Thus, the case (ii) does not occur.
Set $M=\max_{1<j<n}|a_j-a_1|.$
Then we have $e^{\beta_0}M\le \delta_0=|z_0-a_1|.$
If $a_1=0,$ then $M=\tilde\rho_n$ and $\rho_n<e^{\beta_0}\tilde\rho_n\le |z_0|,$
which implies $z_0\in E_n,$ a contradiction.
If $a_1\ne0,$ then $a_k=0$ for some $1<k<n.$
In particular, $|a_1|=|a_k-a_1|\le M.$
Therefore, $e^{\beta_0}M\le |z_0|+|a_1|\le |z_0|+M,$ which yields
$M\le |z_0|/(e^{\beta_0}-1).$
On the other hand, $|a_j|\le |a_j-a_1|+|a_1|\le 2M$ for $1<j<n.$
Therefore, $\rho_n/e=\tilde\rho_n \le 2M\le 2|z_0|/(e^{\beta_0}-1)<|z_0|/e,$
where we used the assumption $\beta_0>\log(1+2e).$
Thus, we obtain $|z_0|>\rho_ne,$ which implies $z_0\in E_n.$
This contradicts again the assumption $z_0\in W.$
Hence, the conclusion follows.
\end{pf}

The next result gives some inequalities for the numbers $\rho_j$ in terms of $\mQ(A).$

\begin{prop}\label{prop:comp}
Let $A$ be a finite subset of $\sphere$ with $0,1,\infty\in A.$
Then,
$$
\frac12\exp\left(\frac12\mQ(A)+2\right)
\le\frac{\rho_n}{\rho_{\min}}
\le2\exp\Big(\mQ(A)+2\Big),
$$
$$
\frac12 \exp\Big(-\mQ(A)-1\Big)\le \rho_{\min}\le e\inv
\aand
e\le\rho_{n}\le \exp\Big(\mQ(A)+1\Big).
$$
\end{prop}

\begin{pf}
Recall that $A=\{a_1,\dots,a_n\}$ with $0,1\in A$ and $a_n=\infty.$
We also recall that $\tilde\rho_j=e\rho_j=\min_{k\ne j}|a_k-a_j|$
for $j=1,\dots,n-1$ and $\tilde\rho_n=\rho_n/e=\max_{j\ne n}|a_j|\ge1.$
We make a preliminary observation.
For a triple $a_j,a_k,a_l$ of distinct finite points in $A,$ we have
$$
\frac{|a_k-a_j|}{|a_l-a_j|}\le \frac{2\tilde\rho_n}{\tilde\rho_j}
\le \frac{2\rho_n}{e^2\rho_{\min}}.
$$
Put $Q=\mQ(A).$
By the definition of $\mQ(A),$ there is a quadruple $a_j,a_k,a_l,a_m$
such that $e^Q=|\CR(a_j,a_k,a_l,a_m)|.$
If none of them is $\infty,$ we estimate $Q$ by
$$
e^Q=\frac{|a_l-a_j|}{|a_k-a_j|}\cdot\frac{|a_k-a_m|}{|a_l-a_m|}
\le \left(\frac{2\rho_n}{e^2\rho_{\min}}\right)^2,
$$
which implies the first inequality.
Otherwise, we have the better estimate $e^Q\le 2a_n/e^2\rho_{\min}.$

Next we show the second inequality.
We may assume that $e\rho_{\min}=|a_2-a_1|$ and $\rho_n/e=|a_3|.$
Then,
$$
\frac{\rho_n}{e^2\rho_{\min}}=\frac{|a_3|}{|a_2-a_1|}
\le|\CR(a_1,a_2,a_3,\infty)|+|\CR(a_1,a_2,0,\infty)|\le 2e^Q,
$$
which leads to the second inequality.
Since $\rho_n\ge e,$ the second inequality yields the third one.
Finally, we have
$$
e\inv\rho_n=|a_3|=|\CR(a_3,\infty,0,1)|\le e^Q
$$
which completes the proof.
\end{pf}

\begin{pf}[Proof of Theorem \ref{thm:q}]
We will show that $K_1h_X(z_1,z_2)\le \qt_X(z_1,z_2)\le K_2h_X(z_1,z_2)$
for $z_1,z_2\in W,$ where $W$ is as in  Theorem \ref{thm:general} and $K_1$ and $K_2$ are constants depending only on $Q(A).$
Then Theorem \ref{thm:general} will provide the required assertion.
Since $\lambda_X(z)\le 1/\delta_X(z)\le 1/\deltat_X(z),$ it is clear that
the inequality $h_X(z_1,z_2)\le\qt_X(z_1,z_2)$ holds for $z_1,z_2\in X.$
In order to show the other inequality, we take a shortest hyperbolic
geodesic $\alpha$ joining $z_1$ and $z_2$ in $X$ for $z_1,z_2\in W.$
By Lemmas \ref{lem:BP} and \ref{lem:beta}, we have the estimate
$$
\lambda_X(z)\ge \frac{K}{\delta_X(z)},\quad z\in W,
$$
where $K=K(Q(A))$ is a constant depending only on $Q(A).$
More concretely, one can choose
$$K=\frac{C_3}{C_1+\max\{Q(A)/2+c,\log(1+2e)\}}\,.$$
Hence,
$$
h_X(z_1,z_2)
=\int_\alpha\lambda(z)|dz|
\ge K\int_\alpha\frac{|dz|}{\deltat_X(z)}
\ge K\qt_X(z_1,z_2)
$$
if $\alpha$ is contained in $W.$
Otherwise, we take a connected component $\alpha_1$ of $\alpha\cap E_j.$
Then the endpoints $\zeta_1,\zeta_2$ of $\alpha_1$ lie on the circle $\partial E_j.$
Choose $a_k$ so that $\tilde\rho_j=|a_k-a_j|~ (\tilde\rho_n=|a_k|$ when $j=n$).
Then the function $w=g(z)=(z-a_j)/(a_k-a_j)$ maps $X$ into $\C_{0,1}$
(let $g(z)=a_k/z$ when $j=n$).
By Lemma \ref{lem:Df} with $f=\id$ and \eqref{eq:qt}, we obtain
\begin{align*}
h_X(\zeta_1,\zeta_2)&\ge h_{\C_{0,1}}(g(\zeta_1),g(\zeta_2))
\ge \frac1{M_0} D^\id(g(\zeta_1),g(\zeta_2))=\frac{\theta}{M_0} \\
&=\frac1{M_0} q(g(\zeta_1),g(\zeta_2))
=\frac1{M_0}q_X(\zeta_1,\zeta_2)
\ge\frac1{2M_0}\qt_X(\zeta_1,\zeta_2),
\end{align*}
where $\theta=|\arg(\zeta_1-a_j)/(\zeta_2-a_j)|\in[0,\pi].$
By using the triangle inequality for $\qt,$ we now obtain
$\qt_X(z_1,z_2)\le K_2 h_X(z_1,z_2)$ with $K_2=\max\{1/K, 2M_0\}.$
\end{pf}

\def\cprime{$'$} \def\cprime{$'$} \def\cprime{$'$}
\providecommand{\bysame}{\leavevmode\hbox to3em{\hrulefill}\thinspace}
\providecommand{\MR}{\relax\ifhmode\unskip\space\fi MR }
\providecommand{\MRhref}[2]{%
  \href{http://www.ams.org/mathscinet-getitem?mr=#1}{#2}
}
\providecommand{\href}[2]{#2}

\end{document}